\documentclass{article}
\usepackage{amssymb,amsmath,stmaryrd,array,graphicx,mathtools,xcolor,academicons}
\definecolor{orcidlogocol}{HTML}{A6CE39}
\title{The price of mathematical scepticism}
 \author{Paul Blain Levy, University of
   Birmingham\footnote{Orcid: 
     {orcid.org/0000-0003-0864-1876} \ Email: \href{mailto:P.B.Levy@bham.ac.uk}{P.B.Levy@bham.ac.uk}}}
\usepackage[colorlinks,breaklinks,citecolor=blue,urlcolor=blue]{hyperref}
\usepackage[numbers]{natbib}
\usepackage{url}
\usepackage[margin=4.1cm]{geometry}
\usepackage[matrix,arrow,curve,rotate]{xy}
\newenvironment{quoti}{\begin{quote}}{\end{quote}}

\newcommand{\distint}[1]{|\!| {#1} |\!|}
\newcommand{\pathree}{Z$_{3}$}
\newcommand{\clof}[1]{\mathsf{Class}({#1})}
\newcommand{\goognat}{\nats_{\mathsf{G}}}
\newcommand{\patwo}{Z$_{2}$}

\newcommand{\Ord}{\mathsf{Ord}}

\input{diagxy}

\newcommand{\pset}{\mathcal{P}}

\newcommand{\subsch}[1]{}

\newcommand{\nats}{\mathbb{N}}

\newcommand{\setbr}[1]{\{{#1}\}}
\newcommand{\tuple}[1]{\langle {#1} \rangle}

\newcommand{\smin}{\!\in\!}

\newdimen\proofrulebreadth \proofrulebreadth=.05em
\newdimen\proofdotseparation \proofdotseparation=1.25ex
\newdimen\proofrulebaseline \proofrulebaseline=2ex
\newcount\proofdotnumber \proofdotnumber=3
\let\then\relax
\def\hfi{\hskip0pt plus.0001fil}
\mathchardef\squigto="3A3B
\newif\ifinsideprooftree\insideprooftreefalse
\newif\ifonleftofproofrule\onleftofproofrulefalse
\newif\ifproofdots\proofdotsfalse
\newif\ifdoubleproof\doubleprooffalse
\let\wereinproofbit\relax
\newdimen\shortenproofleft
\newdimen\shortenproofright
\newdimen\proofbelowshift
\newbox\proofabove
\newbox\proofbelow
\newbox\proofrulename
\def\shiftproofbelow{\let\next\relax\afterassignment\setshiftproofbelow\dimen0 }
\def\shiftproofbelowneg{\def\next{\multiply\dimen0 by-1 }%
\afterassignment\setshiftproofbelow\dimen0 }
\def\setshiftproofbelow{\next\proofbelowshift=\dimen0 }
\def\setproofrulebreadth{\proofrulebreadth}

\def\prooftree{
\ifnum  \lastpenalty=1
\then   \unpenalty
\else   \onleftofproofrulefalse
\fi
\ifonleftofproofrule
\else   \ifinsideprooftree
        \then   \hskip.5em plus1fil
        \fi
\fi
\bgroup
\setbox\proofbelow=\hbox{}\setbox\proofrulename=\hbox{}%
\let\justifies\proofover\let\leadsto\proofoverdots\let\Justifies\proofoverdbl
\let\using\proofusing\let\[\prooftree
\ifinsideprooftree\let\]\endprooftree\fi
\proofdotsfalse\doubleprooffalse
\let\thickness\setproofrulebreadth
\let\shiftright\shiftproofbelow \let\shift\shiftproofbelow
\let\shiftleft\shiftproofbelowneg
\let\ifwasinsideprooftree\ifinsideprooftree
\insideprooftreetrue
\setbox\proofabove=\hbox\bgroup$\displaystyle 
\let\wereinproofbit\prooftree
\shortenproofleft=0pt \shortenproofright=0pt \proofbelowshift=0pt
\onleftofproofruletrue\penalty1
}

\def\eproofbit{
\ifx    \wereinproofbit\prooftree
\then   \ifcase \lastpenalty
        \then   \shortenproofright=0pt  
        \or     \unpenalty\hfil         
        \fi
\fi
\global\dimen0=\shortenproofleft
\global\dimen1=\shortenproofright
\global\dimen2=\proofrulebreadth
\global\dimen3=\proofbelowshift
\global\dimen4=\proofdotseparation
\global\count255=\proofdotnumber
$\egroup  
\shortenproofleft=\dimen0
\shortenproofright=\dimen1
\proofrulebreadth=\dimen2
\proofbelowshift=\dimen3
\proofdotseparation=\dimen4
\proofdotnumber=\count255
}

\def\proofover{
\eproofbit 
\setbox\proofbelow=\hbox\bgroup 
\let\wereinproofbit\proofover
$\displaystyle
}%
\def\proofoverdbl{
\eproofbit 
\doubleprooftrue
\setbox\proofbelow=\hbox\bgroup 
\let\wereinproofbit\proofoverdbl
$\displaystyle
}%
\def\proofoverdots{
\eproofbit 
\proofdotstrue
\setbox\proofbelow=\hbox\bgroup 
\let\wereinproofbit\proofoverdo
ts
$\displaystyle
}%
\def\proofusing{
\eproofbit 
\setbox\proofrulename=\hbox\bgroup 
\let\wereinproofbit\proofusing
\kern0.3em$
}

\def\endprooftree{
\eproofbit 
  \dimen5 =0pt
\dimen0=\wd\proofabove \advance\dimen0-\shortenproofleft
\advance\dimen0-\shortenproofright
\dimen1=.5\dimen0 \advance\dimen1-.5\wd\proofbelow
\dimen4=\dimen1
\advance\dimen1\proofbelowshift \advance\dimen4-\proofbelowshift
\ifdim  \dimen1<0pt
\then   \advance\shortenproofleft\dimen1
        \advance\dimen0-\dimen1
        \dimen1=0pt
        \ifdim  \shortenproofleft<0pt
        \then   \setbox\proofabove=\hbox{%
                        \kern-\shortenproofleft\unhbox\proofabove}%
                \shortenproofleft=0pt
        \fi
\fi
\ifdim  \dimen4<0pt
\then   \advance\shortenproofright\dimen4
        \advance\dimen0-\dimen4
        \dimen4=0pt
\fi
\ifdim  \shortenproofright<\wd\proofrulename
\then   \shortenproofright=\wd\proofrulename
\fi
\dimen2=\shortenproofleft \advance\dimen2 by\dimen1
\dimen3=\shortenproofright\advance\dimen3 by\dimen4
\ifproofdots
\then
        \dimen6=\shortenproofleft \advance\dimen6 .5\dimen0
        \setbox1=\vbox to\proofdotseparation{\vss\hbox{$\cdot$}\vss}%
        \setbox0=\hbox{%
                \advance\dimen6-.5\wd1
                \kern\dimen6
                $\vcenter to\proofdotnumber\proofdotseparation
                        {\leaders\box1\vfill}$%
                \unhbox\proofrulename}%
\else   \dimen6=\fontdimen22\the\textfont2 
        \dimen7=\dimen6
        \advance\dimen6by.5\proofrulebreadth
        \advance\dimen7by-.5\proofrulebreadth
        \setbox0=\hbox{%
                \kern\shortenproofleft
                \ifdoubleproof
                \then   \hbox to\dimen0{%
                        $\mathsurround0pt\mathord=\mkern-6mu%
                        \cleaders\hbox{$\mkern-2mu=\mkern-2mu$}\hfill
                        \mkern-6mu\mathord=$}%
                \else   \vrule height\dimen6 depth-\dimen7 width\dimen0
                \fi
                \unhbox\proofrulename}%
        \ht0=\dimen6 \dp0=-\dimen7
\fi
\let\doll\relax
\ifwasinsideprooftree
\then   \let\VBOX\vbox
\else   \ifmmode\else$\let\doll=$\fi
        \let\VBOX\vcenter
\fi
\VBOX   {\baselineskip\proofrulebaseline \lineskip.2ex
        \expandafter\lineskiplimit\ifproofdots0ex\else-0.6ex\fi
        \hbox   spread\dimen5   {\hfi\unhbox\proofabove\hfi}%
        \hbox{\box0}%
        \hbox   {\kern\dimen2 \box\proofbelow}}\doll%
\global\dimen2=\dimen2
\global\dimen3=\dimen3
\egroup 
\ifonleftofproofrule
\then   \shortenproofleft=\dimen2
\fi
\shortenproofright=\dimen3
 \onleftofproofrulefalse
\ifinsideprooftree
\then   \hskip.5em plus 1fil \penalty2
\fi
}

\begin{document}
\bibliographystyle{plain}
\maketitle

\begin{abstract}
 This paper argues that, insofar as we doubt the bivalence of the
 Continuum Hypothesis or the truth of the Axiom of Choice, we should
 also doubt the consistency of third-order arithmetic, both the 
 classical and intuitionistic versions.

 Underlying this argument is the following philosophical view.
 Mathematical belief springs from certain intuitions, each of which
 can be either accepted or doubted in its entirety, but not
 half-accepted.    Therefore, our beliefs about reality, bivalence,
 choice and consistency should all be aligned.
\end{abstract}


\section{Introduction}
\subsection{Theory
  vs reality}

According to a widely held ``classical'' view of mathematics, the 
Continuum Hypothesis (CH) is \emph{bivalent}, i.e.\ either
objectively true or objectively false, even if it is absolutely unknowable which of these is the case.  Furthermore, the Axiom of
Choice (AC) is true, and therefore also the weaker principle known as Dependent Choice (DC).

There are many other views, however.  Here are two that you may have encountered:
\begin{itemize}
\item {``There is no canonical universe of mathematical reality, but
    rather many universes of equal status.  All of them satisfy the
    ZFC axioms, but CH holds in some of them and fails in others.''}  
\item {``AC is unacceptable because it leads to the Banach--Tarski theorem.  Therefore ZF$+$DC should be adopted as a foundational theory.''} 
\end{itemize}
Each favours a strong foundational theory (at least ZF), yet at the same time is sceptical of the classical conception.

In this article, I shall present an argument that we cannot ``have our cake and eat it'' in this way.   Although scepticism is legitimate, it comes at a price.  

Before this is spelt out, we need some technical preliminaries.

Firstly, let us note that both CH and the Banach--Tarski theorem are \emph{third-order arithmetical} sentences, meaning that---with suitable coding---each quantifier ranges over $\nats$ or $\pset\nats$ or $\pset\pset\nats$, but nothing more complex.  Accordingly, to discuss these sentences, we need not consider advanced theories such as ZF.  Let us merely consider \pathree, the theory of third-order arithmetic.  (Technically, it is a 3-sorted first-order theory, with Extensionality axioms and unrestricted Comprehension and Induction schemes.)  

Secondly, say that a relation $R$ from a set $A$ to a set $B$ is \emph{entire} when, for all $a \smin A$, there is  $b\smin B$ such that $R(a,b)$. 
Then AC and DC are stated as follows.
\begin{description}
\item[AC] For any sets $A$ and $B$, and any entire relation $R$ from $A$ to $B$, there is a function $f \in B^A$ such that, for all $a \smin A$, we have $R(a,f(a))$.
\item[DC] For any set  $B$, and any $b \smin B$ and entire endorelation $R$ on $B$ (i.e.\ relation from $B$ to itself), there is a sequence $(x_n)_{n \in \nats} \in B^{\nats}$ such that $x_0 = b$ and, for all $n \smin \nats$, we have $R(x_n,x_{n+1})$.
\end{description}
The contention of this article is that, insofar as we doubt CH bivalence or AC, we should also doubt the consistency of \pathree{}.  Likewise, doubting DC leads to doubt in the consistency of \patwo{}, the theory of second-order arithmetic.

\paragraph*{Note} In some fields of mathematics, such as topos theory,
it is common to avoid using AC and other classical principles, in
order to gain information about interesting models where these principles fail.  (See~\cite{FiorePittsSteenkamp:quotinduct} for a recent example that actually relies on AC being true in reality.)  Since this practice is not motivated by scepticism, it is philosophically uncontroversial and does not bear on our discussion.

\subsection{Structure of paper}

We proceed as follows.  We begin (Section~\ref{sect:justbelief}) with
a general discussion of belief and doubt, and I set out the paper's fundamental
principles---which are open to dispute, of course.  We then classify
mathematical beliefs (Section~\ref{sect:biv}), based on the bivalence
of different kinds of sentence.  We consider the intuitions that may
give rise to each of the positions (Section~\ref{sect:intu}) and
discuss which of them are reliable (Section~\ref{sect:reli}).  I then
argue that the various ways of answering this question lead to the
claimed consequences (Section~\ref{sect:cons}).  I also consider
intuitionistic theories (Section~\ref{sect:inttheories}), and critique the view that ``reality is indeterminate''
(Section~\ref{sect:multi}).  The issue of encoding is addressed in 
Section~\ref{sect:platnotess}.

We then review some literature, looking at positions similar to the
one I am proposing    
(Section~\ref{sect:similar}) and ones that conflict with it 
(Section~\ref{sect:chall}).

Section~\ref{sect:conc} concludes, with a mention of future work.


\section{Principles of justified belief} \label{sect:justbelief}

The words ``doubt'' and ``scepticism'' have various shades of meaning
in English.  In this article, they refer to a lack of belief in
$X$, not to a belief that $X$ is false.  So please do not interpret me
as saying that CH bivalence or AC sceptics should believe \pathree{} to be inconsistent; they should not.  

The following examples illustrate our basic principles concerning
belief and doubt.  To avoid irrelevant infinity issues, let $\goognat$ be the set of \emph{Googolplex-bounded numbers}, i.e.\ natural numbers less than $10^{10^{100}}$.
\begin{enumerate}
\item Consider the statement ``Cleopatra ate an even number of grapes''.    
  There is no evidence for or against.  A person who believes this may
  happen to be right, but their belief is nevertheless arbitrary and
  unjustified.  So the correct position is to doubt it.
\item  The \emph{Googolplex Goldbach} conjecture says that every even
  Googolplex-bounded number other than 0 and 2 is a sum of two primes.
  We do not know whether this is true, so we doubt it.

  What would cause us to believe it?  Either a proof, or intuition, or
  a combination of the two.  These are (we shall suppose) the only acceptable grounds for belief.  Furthermore, appeals to intuition raise the tricky question of which intuitions are reliable.
  
We might be tempted towards belief by the 
fact~\cite{OliveiraeSilvaHerzogPardi:empverifgc} that the Goldbach
property has been verified up to $4 \times
10^{18}$.  And it is
 {especially} tempting to believe that the property holds for  the least
 even number that has not yet been checked.\footnote{As of 16 January 
   2022, this number is
   $4.01 \times 10^{18} + 4$.  (Personal communication from Tom\'{a}s Oliveira e Silva.)}  But even this
 proposition---call it \emph{Liminal Goldbach}---might be false for
 all we
 know.  So we doubt it.

This illustrates a general principle: inductive
evidence, however strong, is not sufficient grounds for belief.  Mathematicians
throughout the ages have largely agreed on this point.

(A possible objection: it is common mathematical practice to trust a
proof checked by a computer or another person, and this relies on inductive inference.
 We are ignoring such issues.)






\item Consistency statements are not essentially different from
  statements about prime numbers.  Recall that a theory $T$ is
  \emph{consistent} when $\mathsf{False}$ has no $T$-proof.
  Likewise, let us say that $T$ is \emph{Googolplex
    consistent} when $\mathsf{False}$ has no $T$-proof whose
  length is Googolplex-bounded.  I shall assume without further
  comment that proof length has been precisely defined for
  each of our theories. 

  Consider the statement ``\pathree{} is Googolplex consistent''.
  As before, our default position is to doubt it, and only proof or intuition will give us adequate grounds to believe it.  G\"{o}del's second incompleteness theorem and similar results do not justify relaxing this policy.

  One sometimes hears the following argument for
  consistency: ``Many clever people have used this theory and studied
  its foundations for years, and found no contradiction.''  Since this
  is an inductive inference, it is not sufficient grounds for
  belief.\footnote{As Hamkins~\cite{Hamkins:smart} points out: the
    negation of Fermat's Last Theorem turned out to be inconsistent,
    even though,  before Wiles, many clever people had looked seriously and
           been unable to refute it.}
\end{enumerate}
To summarize:
\begin{itemize}
\item For any statement, our default position is doubt.
\item Only proof and/or intuition will move us to a state of belief.
\item We need to decide which intuitions are reliable.
\item Inductive inference is not accepted.
\item These principles apply, in particular, to consistency statements.
\end{itemize}



\section{The bivalence questionnaire} \label{sect:biv}

Our investigation of mathematical reality begins with a list of
sentences whose truth value is unknown.  The details of each sentence
are not so important, but please pay
attention to the \emph{logical form}.  The list is as follows.
\begin{description}
\item[Physical sentences] which concern the physical universe.
  \begin{itemize}
  \item The \emph{Cleopatra Hypothesis}: Cleopatra ate an even number
    of grapes.
\end{itemize}
\item[Computational sentences] where quantifier ranges are finite.
  \begin{itemize}
  \item The \emph{Googolplex Goldbach conjecture}: Every even
    Googolplex-bounded number other than 0 and 2 is a sum of two
    primes.
\end{itemize}
\item[Arithmetical sentences] where quantifiers range over all natural numbers.
  \begin{itemize}  \item 
    The \emph{Goldbach conjecture}~\cite{wikipedia:goldbach}: Every
    even natural number other than 0 and 2 is a sum of two
    primes. This has the form $\forall n \smin \nats. \,\phi(n)$, where
    $\phi$ is computational.
  \item 
    The \emph{twin prime conjecture}~\cite{wikipedia:twinprime}: There
    are infinitely many $n \smin \nats$ such that both $n$ and $n+2$ are
    prime.  This has the form $\forall n \smin \nats.\,\exists m \smin
    \nats. \phi(m,n)$, where $\phi$ is computational.
  \end{itemize}
\item[Second-order arithmetical sentences] where quantifiers range over all sets of natural numbers.
  \begin{itemize}
  \item The \emph{Littlewood conjecture}~\cite{wikipedia:littlewood}:
    For any real numbers $\alpha$ and $\beta$, we have  \linebreak $\liminf_{n
      \rightarrow \infty} n \distint{n \alpha} \distint{n \beta} = 0$,
    where $\distint{\,}$ is the distance to the nearest integer.  This
    has the form $\forall x \smin \pset \nats.\,\phi(x)$,
    where $\phi$ is arithmetical, because a real number can be encoded
    as a set of natural numbers, and a pair of subsets of $\nats$ can
    be encoded as a single subset.
  \item The \emph{Toeplitz conjecture}~\cite{wikipedia:inscribedsq}:
    Every simple closed curve contains all four vertices of some
    square.  This has the form $\forall x \smin
    \pset\nats.\,\exists y\smin \pset \nats.\,\phi(x,y)$, where $\phi$
    is arithmetical, because a continuous function on $\mathbb{R}$
    can be encoded as a continuous function on $\mathbb{Q}$.
  \end{itemize}
 \item[Third-order arithmetical sentences] where quantifiers range over all sets of sets of natural numbers.
   \begin{itemize}
   \item The \emph{Continuum Hypothesis}~\cite{wikipedia:ch}:  There
     is a bijection from $\aleph_1$ to $2^{\nats}$.  This has the form $\exists x \smin
     \pset\pset\nats.\,\phi(x)$, where $\phi$ is second-order
     arithmetical, because an element of $\aleph_1$ can
     be encoded (non-uniquely) as a well-ordered subset of $\nats$.
   \item The \emph{Suslin Hypothesis}~\cite{wikipedia:suslin,Glazer:pitwotwo}: The
     tree $\setbr{0,1}^{<\omega_1}$ has no subtree in which every
     chain and every antichain is countable.  This has the form
     $\neg\exists x \smin \pset\pset\nats. \forall y \smin
     \pset\pset\nats.\,\phi(x,y)$, where $\phi$ is second-order
     arithmetical, because an element of the tree can
     be encoded (non-uniquely) as a tuple $\tuple{X,\prec,Y}$, where
     $\tuple{X,\prec}$ is a well-ordered subset of $\nats$, and $Y
     \subseteq X$ indicates which elements of $X$ are mapped to 1.
     \end{itemize}
\item[Unrestricted set-theoretic sentences] where quantifiers range
  over all sets or all ordinals.
  \begin{itemize}
  \item The \emph{Generalized Continuum
      Hypothesis}~\cite{wikipedia:ch}: Every infinite cardinal
    $\kappa$ satisfies $\kappa^+ = 2^{\kappa}$.  This has the form
    $\forall \alpha \smin \Ord.\,\phi(\alpha)$, where $\phi$ uses only
    restricted quantifiers,\footnote{The sense is more liberal than
      that of the limited
      ZF syntax, and allows $\pset$.} since cardinals can be
    encoded as initial ordinals.
  \item The \emph{Eventually Generalized Continuum Hypothesis}:  There
    is an infinite cardinal $\lambda$ such that every cardinal $\kappa
    \geqslant \lambda$ satisfies  $\kappa^+ = 2^{\kappa}$.   This has
    the form  $\exists \alpha \smin \Ord.\,\forall \beta \smin \Ord.\,\phi(\alpha,\beta)$, where $\phi$ uses only
    restricted quantifiers.
  \end{itemize}
\item[Class-theoretic sentences] where quantifiers range over all classes.
  \begin{itemize}
  \item The \emph{Club-Failure Hypothesis} considered by Schlutzenberg~\cite{Schlutzenberg:pioneone}:
    Every club class of infinite cardinals has a member whose
    successor cardinal $\kappa$ is a GCH failure, i.e.\ satisfies $\kappa^+ < 2^{\kappa}$.    Writing
    $\clof{\Ord}$ for the collection of classes of ordinals, this has
    the form $\forall X \smin \clof{\Ord}.\,\phi(X)$, where $\phi$ is set-theoretic. 
  \item The \emph{Ord--Suslin Hypothesis} considered by Hamkins and Switzer~\cite{Switzer:pionetwo}: The tree
    $\setbr{0,1}^{<\Ord}$ has no subtree in which every chain and
    every antichain is a set.  This has the form $\neg \exists X \smin
    \clof{\Ord}.\,\forall Y \smin \clof{\Ord}.\,\phi(X,Y)$, where $\phi$
    is set-theoretic.
  \end{itemize}
\end{description}
For each kind of sentence other than physical and computational, I
have given two examples.  The first has just one quantifier
($\forall$ or $\exists$) of the specified kind, and the second has two
($\forall \exists$ or $\exists \forall$).  Such sentence forms are often described in the
terminology of logical complexity: the Goldbach conjecture is
$\Pi^0_1$, the twin prime conjecture $\Pi^0_2$, the Littlewood
conjecture $\Pi^1_1$ and so forth.  Note that the Goldbach conjecture
(like every $\Pi^0_1$ sentence) is ``falsifiable'', meaning that
someone who asserts it runs the risk of being hit by a
counterexample.    

Now I am going to interrogate you.  Assume pessimistically that these
sentences cannot be proved or refuted \emph{in any convincing way} within the
lifetime of the universe.  
  Under this assumption (which may be correct for all we know), which
  of these sentences do you consider to
  be bivalent?  In other words, do you think that---despite our hopeless
  ignorance---there is a fact of the matter 
  whether Cleopatra ate an even number of grapes?  Whether every even 
  Googolplex-bounded number other than 0 and 2 is a sum of two primes?
  And so forth.

  Let me stress that ``bivalence ambivalence'' is allowed and even 
  encouraged.  I dare say that few people would be certain of
  their answer for all the sentences.  

Questionnaires of this kind have often
appeared~\cite{Gaifman:ontology,Velleman:CHpost}.  They provide a crude but useful device to measure a
person's belief in objective reality, a belief known as
``realism'' or ``platonism''.  (These words will be used 
   interchangeably.)    

It is time to name various philosophical positions. 
\begin{itemize}
\item An \emph{ultrafinitist}~\cite{vanDantzig:finitenumber,Parikh:existfeas} doubts that
  computational sentences are bivalent. 
\item A \emph{finitist} accepts this, but doubts that arithmetical
  sentences are bivalent.  
\item A \emph{countabilist} accepts this, but doubts that second-order
  arithmetical sentences are bivalent. 
\item A \emph{sequentialist} accepts this and also DC, but doubts that
  third-order arithmetical sentences are bivalent.
\item A \emph{particularist} accepts this and also AC, but doubts that
  sentences that quantify over all sets or ordinals are bivalent.
\item A \emph{totalist} accepts all the above, but doubts that sentences that quantify over all classes are bivalent.  
\end{itemize}
The above taxonomy immediately raises questions. Is this all just a choice
between various coloured pills?  Why not have an option for someone
who accepts that 17th order arithmetical sentences are bivalent but
not 18th order?  Or for someone who accepts that $\Pi^0_{52}$
sentences are bivalent but not $\Pi^0_{53}$ ones?  

To answer these questions, recall that belief should not be arbitrary.
Furthermore (according to our principles), it should be justified by proof or intuition.  So we cannot be mere ``truth value realists'', believing for no reason that sentences of a certain kind are bivalent.  
 What, then, are the intuitions that support the various positions?

\section{Intuitions of mathematical reality} \label{sect:intu}



I now present five intuitions that I experience, and hopefully you do too.  They are little people inside our head, and each of them is going to speak.  For the moment, just listen to them.  We postpone the question of whether they are reliable. 
\begin{description}
\item[Googolplex] 
  {``I perceive the notion of Googolplex-bounded number.   Since this is a clearly defined notion,  quantification over the set $\goognat$ yields an objective truth value.''}
\item[Arbitrary Natural Number]  {``I perceive the notion of a natural number, given by zero and successor and nothing more.   This is a clearly defined notion, as restrictive as possible.  So quantification over the set $\nats$ yields an objective truth value.''}
\item[Arbitrary Sequence] {``Given a set $B$, I perceive the notion of a sequence $(x_n)_{n \in \nats}$ in $B$, which consists of successive arbitrary choices of an element of $B$.   This is a clearly defined notion, as liberal as possible.  So quantification over the set $B^{\nats}$ yields an objective truth value.  Since a sequence consists of successive arbitrary choices, DC holds.''}   
\item[Arbitrary Function]  
  {``Given sets $A$ and $B$, I perceive the notion of a function $f$ from $A$ to $B$, which consists of independent arbitrary choices $f(a) \in B$, one for each $a \in A$.  This is a clearly defined notion, as liberal as possible.  So quantification over the set $B^A$ yields an objective truth value.  Since a function consists of  independent arbitrary choices, AC holds.''}
\item[Arbitrary Ordinal]
  {``I perceive the notion of an ordinal.  This is a clearly defined notion, as liberal as possible.  So quantification over the class $\Ord$ yields an objective truth value.''}
\end{description}
Let me stress that AC is integral to the intuition of a function as
consisting of independent arbitrary choices; it is not a separate
intuition.  Likewise, DC is integral to the intuition of a sequence as
consisting of successive arbitrary choices.


Ultrafinitists accept none of these intuitions.  Finitists accept the
first one, countabilists the first two, sequentialists the  first
three, particularists the first four, and totalists all
five.  
  To see the link to \patwo{} and \pathree{}, note that for any set $A$ we have a bijection $\pset A \cong \setbr{0,1}^A$ that represents each subset $C$ of $A$ by its characteristic function.\footnote{Constructivists would say that I am assuming $C$ to be a ``decidable'' subset, meaning that every $a \smin A$ is either a member or a non-member of $C$. 
  Since our concern at this point is bivalence and classical theories, that is not an issue.}
 So we have bijections $\pset \nats \cong \setbr{0,1}^{\nats}$ and $\pset\pset
 \nats \cong \setbr{0,1}^{\setbr{0,1}^{\nats}}$.    

The taxonomy I have given is crude, and ignores many important distinctions.  For example, it
is usual~\cite{Wang:eightyyearsfoundstud} to distinguish finitists
from 
 {constructivists} (also called intuitionists).  The latter believe
  in \emph{constructions} on natural numbers, and in 
higher-order constructions; but such notions are beyond the scope of this
paper.   I have also avoided the
question of what ultrafinitists \emph{do} believe, but presumably
 calculations are acceptable.\footnote{Friedman~\cite{Friedman:philprobmath} recalls  meeting the ultrafinitist Esenin-Volpin: ``I then proceeded to
  start with $2^1$ and asked him whether this is `real' or something
  to that effect. He virtually immediately said yes. Then I asked
  about $2^2$, and he again said yes, but with a perceptible
  delay. Then $2^3$, and yes, but with more delay. This continued for
  a couple of more times, till it was obvious how he was handling this
  objection. [\ldots] There is no way that I could get very far with this.''}   

Finitism is linked to Primitive
Recursive Arithmetic (PRA), a subsystem of PA~\cite{Tait:finitism}.  
  Countabilism is linked to the ``predicative'' line of work initiated by
  Weyl~\cite{Weyl:daskont}.  Sequentialism
  is linked to descriptive set
  theory~\cite[Section 3.2.3]{Rathjen:scopefeferman}.

\section{Reliability of the intuitions} \label{sect:reli}

We now need to consider which of these intuitions are reliable.  Each one claims to have (limited) access to an objective, ``platonic'' realm, much larger than we can directly apprehend. 
So the truly sceptical answer is that none are reliable.  \emph{``How
  can a human mind have access to an immense platonic realm?  The idea
  is absurd!''}\footnote{Such arguments have been considered by
  Benacerraf~\cite{Benacerraf:mathtruth},
  Field~\cite{field:realmathmod} and others.}  The price of such an attitude is ultrafinitism.

This is a contentious point, because finitists and constructivists
sometimes argue in exactly this way against more credulous positions.
But the anti-platonist argument has nothing to do with infinity
\emph{per se}.  The set $\goognat$ is no more capable of direct
apprehension, by an actual human or computer, than $\nats$.  To rescue
finitism from the charge of being a platonist philosophy, some may say
that Googolplex Goldbach can ``in principle'' be decided by checking,
while others may ``prove'' its bivalence by an induction up to
Googolplex.  But surely each of these defences relies on a prior
belief in the very set $\goognat$ that ultrafinitists
doubt. 

In summary, anyone who believes in the bivalence of Googolplex
Goldbach is a platonist.  
 Welcome to the club!

Leaving aside ultrafinitists, then, we all accept $\goognat$ and are
platonists.  Now we must decide how far to go, and it is not an easy
question.   The first four intuitions have the following profound differences.
\begin{itemize}
\item The set $\goognat$ can in principle be grasped.
\item Each element of $\nats$ can in principle be grasped.
\item An element of $B^{\nats}$ is given by just one choice at a time,
  and each time-point can in principle be grasped.
\item An element of $B^A$ is given by $A$-many choices at the same time.
\end{itemize}
Suppose we accept Arbitrary Natural Number and Arbitrary Sequence.
Shall we accept Arbitrary Function?  Two arguments have been made against it.

Firstly, some have suggested that the independence
results~\cite{Goedel:consistCH,Cohen:settheoryCH,LevySolovay:meascardCH}
provide evidence against CH bivalence.   I do not see why that should
be so, even if the truth value of CH is absolutely unknowable.  Whether Cleopatra ate
an even number of grapes is unknowable, but that is not an argument
against bivalence.\footnote{Gaifman~\cite{Gaifman:ontology} makes
  a similar point: ``I toss a fair coin and, without noting
  the side it landed, I toss it again.  Nobody knows and nobody will be
  able to know the outcome of the first toss.  [\ldots] Yet people
  have no problem in believing that there is a true answer to the
  question `which side did the coin land?'{}''.}

In fact, none of the intuitions claim to have \emph{complete}
knowledge of the entities they perceive.  On the contrary, they
profess extreme ignorance, merely claiming to know the most basic
properties.  Most of our knowledge about $\nats$, for example, comes
from proof, not directly from the Arbitrary Natural Number intuition.  Whatever limits may exist on our proof ability, they do not call into question the reliability of that intuition.

This point also applies in reverse.  If the CH mystery is solved at
some future time, this will not give us a reason to consider the
Arbitrary Function intuition reliable.  The one has nothing to do with the other.

Secondly, some have suggested that the Banach--Tarski theorem provides evidence against AC.  But this criticism is based on geometric intuition, which mathematicians have learnt to distrust.  Furthermore, it has been argued that there are also theorems provable without AC that violate geometric intuition~\cite{Feferman:monsters}. 

Discounting these arguments against the Arbitrary Function intuition, we are still left with the question of whether to accept it.  Personally I find the intuition strong enough to accept, but am not free of ambivalence, and can understand others being more cautious.

Lastly we come to the Arbitrary Ordinal intuition.  It is highly
controversial~\cite{Parsons:setsandclasses,Dummett:fregephilmath,RayoUzquiano:absgen,Florio:unrestricted,Rumfitt:boundarythought} 
because of the Burali-Forti paradox: it claims to perceive a notion of
ordinal that is ``as liberal as possible'' and yet excludes the
order-type of the well-ordered class $\Ord$.  Although it has been
defended~\cite{Boolos:whencecontra,Cartwright:speakingevery}, I
personally am sceptical.  In any case, the issue is beyond the scope
of this article, which is concerned with higher-order arithmetic.

\section{Consequences of belief and doubt} \label{sect:cons}
\subsection{Introduction}

Here is the summary so far.  The starting position was that the only
acceptable grounds for belief are proof and intuition.  I then listed
some intuitions that I experience (and am assuming that there are no
others that would undermine my argument).  The key question was which
of these are reliable, and we have noted various possible answers.
In this section, we consider their consequences.


\subsection{Peano Arithmetic and beyond}

We begin with the Arbitrary Natural Number intuition.   If we accept
it, then we believe in a platonic realm of natural numbers, and the bivalence of every arithmetical statement.  Every PA axiom is true, and every inference rule preserves truth.  So every theorem is true, and PA is consistent.

If, on the other hand, we doubt Arbitrary Natural Number, i.e.\ we are
finitists, then this simple path to PA consistency is blocked.  But
perhaps some other proof will convince us.   
So we turn to the literature.  Gentzen~\cite{Gentzen:consistarith,Chow:consistarith}
proved PA consistency using \emph{induction up to $\varepsilon_0$},
and G\"{o}del's Dialectica
argument~\cite{AvigadFeferman:dialectica,goedel:dial} proved it using
\emph{higher-order constructions}.  Unless we accept one of these
principles, we have to doubt the Googolplex consistency of PA.  

Now consider the theory \patwo{}.  We believe it to be consistent if we accept
Arbitrary Natural Number and Arbitrary Sequence.  On the other hand,
if we stop at Arbitrary Natural Number, i.e.\ we are countabilists, then the simple consistency
argument does not work.  But perhaps some other proof will convince
us.  So we turn to the literature.  Tait~\cite{Tait:nonconst} and Girard~\cite{Girard:phdthesis}
proved it using a predicate on $\nats$ defined ``impredicatively'' via quantification over
$\pset\nats$, but a countabilist would surely not accept such a
definition.  
There is also Spector's proof of \patwo{} consistency, which uses 
\emph{higher-typed bar
  recursion}~\cite{Ferreira:Spectorproof,Spector:barrecursion}.
Unless we accept this principle---which is rather unlikely---we have to doubt the Googolplex consistency of \patwo{}.

Moving on, if we accept Arbitrary
Function, then we conclude that CH is bivalent and \pathree{} is
consistent.  But if we stop at Arbitrary Sequence, i.e.\ we are
sequentialists, then we have to doubt the Googolplex consistency of
\pathree{}.  There is no middle ground.  A formal
consistency proof might convince us, but the ones in the
literature---e.g.~\cite{Prawitz:haupt,Takahashi:proofcutelim}---use a 
predicate on $\pset\nats$ defined ``impredicatively'' via
quantification over $\pset\pset\nats$, which a sequentialist would surely not
accept.

My key point is that, although we can either accept or doubt an intuition, we cannot half-accept.  If we consider an intuition to be unreliable, then we should fully discard it.  For a historical example: once the mathematical community came to view geometrical intuition as unreliable, it was fully discarded, in the sense that appealing to it in a proof was no longer allowed.

Thus, accepting \pathree{} but not AC is not an option, since AC is
asserted by Arbitrary Function.  So if the Banach--Tarski
theorem is anything less than an objectively true statement, then either
Arbitrary Natural Number or Arbitrary Function is an unreliable intuition, and the Googolplex consistency of \pathree{} is in doubt.

Likewise, it cannot be said that $\pset^{73}\nats$ is ``less
credible'' than $\pset \pset \nats$.   This is a statistical way of
thinking, appropriate only when inductive inference is used.   If
$\pset^{73}\nats$ is not real, then either Arbitrary Natural Number or
Arbitrary Function is an unreliable intuition, so the reality of $\pset\pset\nats$ is in doubt.

In the same way, our line of thinking does not allow the
``positivist'' view of Kahrs~\cite{Kahrs:formalistpers}, which accepts
the bivalence of $\Pi^0_1$ sentences  but doubts that of the twin
prime conjecture.  For if Arbitrary Natural Number is an unreliable
intuition, then even the bivalence of the Goldbach conjecture is in
doubt.   Likewise for each pair of sentences in our questionnaire: if
we doubt the bivalence of the second sentence, then we should also
doubt that of the first.

\subsection{Finitism and ultrafinitism} \label{sect:finconst}

Let us revisit the two most sceptical schools, which differ in 
their view of PRA.  For a finitist, each PRA axiom is true and each proof
rule preserves truth.  So each theorem is true and PRA is
consistent.  
But an ultrafinitist cannot accept
this argument.  Presumably they will also be
unconvinced by the formal 
consistency proof, which uses induction up to
$\omega^{\omega}$. So they will doubt the Googolplex consistency of
PRA.

Finally, let us note that  a finitist cannot accept the bivalence of  $\Pi^0_1$ statements
such as ``PA is consistent'' or ``ZFC is consistent''.
And an ultrafinitist cannot even accept the bivalence of computational statements such as
``PRA is Googolplex consistent'' or ``ZFC is Googolplex consistent''.

\section{Intuitionistic theories} \label{sect:inttheories}

Up to this point, the theories we have seen are \emph{classical}, i.e.\ they include the law of
Excluded Middle $\phi \vee \neg \phi$.  Dropping this law from PA gives the
intuitionistic theory  
known as \emph{Heyting arithmetic} (HA).  Dropping it from \patwo{}
gives \emph{intuitionistic second-order arithmetic} (I\patwo{}), and
dropping it from \pathree{} gives \emph{intuitionistic third-order
  arithmetic} (I\pathree{}).   Furthermore, we can obtain 
\emph{intensional} versions of I\patwo{} and I\pathree{} by dropping
Extensionality.  The following results are provable in PRA~\cite[page
170]{TroelstravanDalen:constmathone}.
\begin{itemize}
\item PA and HA are equiconsistent.
\item \patwo{} and intensional I\patwo{} are equiconsistent.
\item \pathree{} and intensional I\pathree{} are equiconsistent.
\end{itemize}
We accordingly ask: is there an intuition---other than the ones on our
list---that would yield the
consistency of these intuitionistic theories (and therefore also the
classical ones)?   In the case of HA, the previously mentioned
notion of higher-order construction may be considered such an intuition.

But for I\patwo{} (or intensional I\patwo{}), the answer seems to be
no.  The problem is that the Comprehension scheme allows 
  quantification ranging over $\pset\nats$, which is impredicative.  It is hard
to see how someone who doubts the Arbitrary Natural Number or Arbitrary
Sequence intuition can justify this.\footnote{This concern does not apply to theories
  without powerset, such as dependent type theory with predicative
universes~\cite{MartinLoef:inttypetheory} or
CZF~\cite{Crosilla:czfizf}.    These theories are weaker in
consistency strength than \patwo{}, and are supported by proof-theoretic
results~\cite{Setzer:upperboundmlttwuniv} and constructive
intuitions.  One can of course debate whether these provide sufficient
justification, but the issue is beyond the
scope of this paper.}

Likewise for I\pathree{} (or intensional I\pathree{}), the answer
seems to be no, because the Comprehension scheme allows 
 quantification ranging over 
 $\pset\pset\nats$, which is impredicative.  It is hard to
see how someone who doubts the Arbitrary Natural Number or Arbitrary
Function intuition can justify this.

Some authors have  
proposed \emph{free topos theory} (the axiomatic theory of an
elementary topos with a natural
numbers object) as a foundation of
mathematics~\cite{Lambek:whatworldmath,LambekScott:reflectcat}.  It
contains I\pathree{}, so---according to our
argument---we should doubt its consistency if we doubt CH bivalence
or AC.\footnote{Cf.\  the discussion in~\cite[Section 1.7.2]{LambekScott:reflectcat}.}







\section{Multiversism} \label{sect:multi}

Let us next examine a particular kind of bivalence scepticism, called
\emph{multiversism}~\cite{AntosFriedmanHonzikTernullo:multivconcset,Hamkins:multiverse}.
It asserts that there are many mathematical universes, all of equal status.  In short, ``reality is indeterminate''.  Supposedly, a non-bivalent sentence is one that holds in one universe and not in another.  Sometimes an analogy is drawn with Euclid's Fifth Axiom. 

I shall now raise two concerns with multiversism.\footnote{A 
multiverse theory can be construed either as a philosophical view of
reality or as a mathematical account of a class of models.  My
comments only concern the former.}

Firstly, there is usually a theory (such as ZFC) that all the universes are supposed to model.  This theory needs to be consistent, or else the multiverse will be a ``nulliverse''.  As discussed above, it is hard to see how the belief in consistency can be justified.

Secondly, it seems that multiversism fails to properly account for bivalence doubt.  Let me give some examples of this.  
\begin{itemize}
\item A finitist's doubt in the bivalence of the Goldbach conjecture 
  stems from a fear that $\nats$ may be unreal, not from a fear
  that it may be indeterminate.  In other words, the finitist does not fear that there may 
  be several versions of $\nats$, with the Goldbach conjecture holding
  in one and failing in another.  For in their view, if the conjecture holds in some 
  ``version of $\nats$'' that is at least a model of Robinson
  arithmetic (say), then it is simply
  true.
\item A countabilist's doubt in the bivalence of the Littlewood conjecture  
  stems from a fear that $\pset\nats$ may be unreal, not from a fear
  that it may be indeterminate.   For in their view, if the conjecture fails 
  in some 
  ``version of $\pset\nats$'' that is at least a collection of subsets
  of $\nats$, then it is simply false. 
\item A sequentialist's doubt in CH bivalence stems from a fear that $\pset\pset\nats$ may be unreal, not from a fear that it may be indeterminate.  
  For in their view, if CH holds in some ``version of $\pset\pset\nats$'' that is at
  least a
  collection of subsets of $\pset\nats$, then it is simply true. 
\item A particularist's doubt in GCH bivalence stems from a fear that $\Ord$  may be unreal, not from a fear that it may be indeterminate.  
  For in their view, if GCH fails in some ``version of $\Ord$'' that is at
  least a collection of ordinals, then it is simply false. 
 \item A totalist's doubt in the bivalence of the Club-Failure Hypothesis stems from a fear that
   $\clof{\Ord}$ may be unreal, not from a fear that it may be
   indeterminate.  For in their view, if the hypothesis fails in some ``version of
   $\clof{\Ord}$'' that is at least a collection of classes of
   ordinals, then it is simply
   false.  
\end{itemize}
Each of these examples concerns a single-quantifier sentence just
beyond the boundary of platonist belief.  Doubt in the bivalence of
such sentences cannot be attributed to a fear of
indeterminacy.\footnote{Cf.\ the discussion of $\Pi^0_1$ sentences in Koellner~\cite{Koellner:truthmathqplural}.}

\section{Platonism is not essentialism} \label{sect:platnotess}

I have now almost finished presenting my argument.  In brief, it says
that---with the possible exception of formal consistency proofs---the only adequate basis for
consistency belief is platonism.  My final task is to correct a
certain 
 misunderstanding of platonism that makes my position    
seem more demanding than it actually is.  The issue (or a version of it) has been raised by   
 Benacerraf~\cite{Benacerraf:numbersnot},
Reynolds~\cite{Reynolds:typesabstraction} and others.

To illustrate the point, here are two injections from $\mathbb{Q}$ to $\nats$.  The \emph{red encoding} sends
0 to 0, and $\frac{m}{n}$ (for coprime $m,n >0$) to $2^m\times 3^n$, and
$-\frac{m}{n}$ (for coprime $m,n >0$) to $2^m\times 3^n \times 5$.   The \emph{yellow encoding} sends
0 to 17, and $\frac{m}{n}$ (for coprime $m,n >0$) to $2^m \times 3^n
\times 5$, and
$-\frac{m}{n}$ (for coprime $m,n >0$) to $2^m \times 3^n$.  The 
choice between them is arbitrary, as there is no
reason to prefer one to the other.  Now consider the
following statements:
\begin{itemize}
\item \label{item:redbel} ``Via the red encoding, any believer in the reality of $\nats$
  must also believe in the reality of $\mathbb{Q}$.''
  \item \label{item:yellbel} ``Via the yellow encoding, any believer in the reality of $\nats$
  must also believe in the reality of $\mathbb{Q}$.''
\end{itemize}
According to  a naive construal of
platonism---let us call it \emph{essentialism}---neither statement is
acceptable, and certainly not both.  That is because ``believing in the
reality of $\mathbb{Q}$'' supposedly requires a metaphysical
commitment that goes 
beyond   
 merely  
noting an arbitrary encoding of $\mathbb{Q}$ into another totality that one believes to be
real.   For example, the rational number $\frac{1}{2}$ cannot be
captured ``in its essence'' by its 
 red encoding (18) or its yellow encoding (90).

Yet surely everyone would accept the above
statements.  That is why, although people sometimes say ``I believe in
the reality of $\nats$ but have doubts about $\mathbb{R}$,'' nobody
ever says ``I believe in the reality of
$\nats$ but have doubts about $\mathbb{Q}$,'' and we would be
astonished to hear this.    
So the construal of platonism as essentialism is evidently incorrect.
  Platonists do not, in fact, think that the rational number 
$\frac{1}{2}$ has some kind of transcendent essence, and  
would be content with either encoding of $\mathbb{Q}$.

In summary, we must take care in construing the phrase ``believing in the
reality of $X$'', where $X$ is a totality such as
$\goognat$,  $\nats$, $\mathbb{Q}$, $\mathbb{C}$, $\pset\pset\nats$ or
$\mathsf{Ord}$.   Mere truth value realism is too little, but
essentialism is too much.  








\section{Similar views} \label{sect:similar}

Let us now
review some literature.  Section~\ref{sect:chall}
will consider various critical views and challenges, but first I shall
point to views that agree with aspects of my argument.



\subsection{Finitism vs ultrafinitism}

In Section~\ref{sect:reli}, I described finitism as a platonist
philosophy, because of its belief in the reality of the 
set $\goognat$.  Here
are some similar views, beginning with van Dantzig~\cite{vanDantzig:finitenumber}:
\begin{quoti}
  Unless one is willing to admit fictitious ``superior minds''
  [\ldots] it is necessary, in the foundations of mathematics like in
  other sciences, to take account of the limited possibilities of the
  human mind and of mechanical devices replacing it.
\end{quoti}
Bernadete~\cite[page 210]{Bernadete:book}:
\begin{quoti}
  This standard [finitist] concept of even the potential infinite is no
less dubious than our standard concept of the actual infinite, when
called to account by [ultrafinitist]
proto-mathematics.
\end{quoti}
Wang~\cite{Wang:eightyyearsfoundstud}:
\begin{quote}
  Finitism [...] is an
  idealization.
\end{quote}
Kreisel~\cite{Kreisel:reviewwitt}:
\begin{quote}
  Finitism is, of course, an idealization.
\end{quote}
Troelstra
and van Dalen~\cite[page 6]{TroelstravanDalen:constmathone}:
\begin{quoti}
  All the constructivist schools described in section 1
[including finitism] contain elements of idealization.
\end{quoti}
Nelson~\cite{Nelson:warnsign}:
\begin{quote}
  Finitism is the last refuge of the platonist.
\end{quote}

\subsection{Doubting consistency}

I can hardly claim that this paper's main contention---lack of platonist belief leads to
consistency doubt---is new.  

Firstly, it is the very basis of Hilbert's programme and subsequent work on formal consistency proof~\cite{RathjenSieg:prooftheory}.  That programme was an attempt to convince reality sceptics (specifically, finitists) that a theory expressing certain platonist beliefs (PA) is in fact consistent.  The underlying assumption is that their lack of platonist belief leaves the sceptics with insufficient basis to believe that the theory is consistent.  If that were not so, why bother to look for a consistency proof?   G\"{o}del's second incompleteness theorem tells us that the search is futile (unless the sceptics accept some proof principle that is outside the theory), but does not change the predicament.

Secondly, here are some similar views, beginning with 
 Dummett~\cite{Dummett:whatismathabout}:
\begin{quoti}
  Why, then, does he [the nominalist Field] believe ZF to be consistent? Most people do,
indeed: but then most people are not nominalists. They believe ZF to
be consistent because they suppose themselves in possession of a
perhaps hazily conceived intuitive model of the theory; but Field can
have no such reason.
\end{quoti}
Parsons~\cite[page 58]{Parsons:maththoughtobj}:
\begin{quoti}
  It is [\ldots] hard to see what grounds other than inductive the
  nominalist can have for believing consistency statements for
  theories having only infinite models to be true.
\end{quoti}
Koellner~\cite{Koellner:fefermanonsetth}:
\begin{quoti}
 {I think that
 the concept [suggested by Feferman] of being clear enough to secure consistency (and what the structure is \emph{supposed} to be like) but not clear enough to secure definiteness is itself inherently
 unclear.}
\end{quoti}
D\v{z}amonja~\cite{KantDzamonja:interview}:
\begin{quoti}
  The ZFC axioms [\ldots] have other models too.  But, somehow,
  believing in the consistency comes back to thinking if there is this
  universe of sets or not.
\end{quoti}
Potter~\cite{Potter:inaccess}:
\begin{quoti}
  {As soon as we accept the image of God constructing set
 theory, and exercising free choice in how He constructed it, we must allow the
 possibility that He chose not to construct it at all.}
\end{quoti}
Lastly, here are some specific consistency doubts.
\begin{itemize}
\item The consistency of \pathree{} was doubted by Silver~\cite{Solovay:Silverattest}.
 \item The consistency of \patwo{} was doubted by Gentzen~\cite{Gentzen:presstate},
Lorenzen~\cite{Lorenzen:constrmathphilprob} and P\'{e}ter~\cite[page
233]{Peter:playingwithinf}.
\item The consistency of PA was doubted by the finitist Goodstein~\cite{GoodsteinLakatos:fom}.
\item The consistency of PRA was doubted by the ultrafinitist Nelson~\cite{Nelson:incons}.
\end{itemize}

\subsection{Axiom of Choice}

In Section~\ref{sect:intu}, I described AC as an integral part of the
Arbitrary Function intuition that underpins belief in the bivalence of
third-order arithmetical sentences.  Here are some similar views,
beginning with Jourdain~\cite{Jourdain:reviewyoung}:
\begin{quoti}
  The multiplicative axiom [AC] is necessary in order to be
    able to say that $x,y,z,\ldots$ (an infinity) has any meaning at
    all.
  \end{quoti}
  Bernays~\cite{Bernays:platonism}:
  \begin{quoti}
    The axiom of choice is an immediate application of the
    quasi-combinatorial concepts in
    question.
  \end{quoti}
Shoenfield~\cite{Shoenfield:axioms}:
  \begin{quoti}
    If we interpret a collection as being an arbitrary division of the
    objects available into members and non-members of the collection,
    it is reasonable to claim that such a collection [representing a
    choice function] exists.
  \end{quoti}
Ferreir\'{o}s~\cite{Ferreiros:arbsetszfc}:
\begin{quoti}
    From the standpoint of a principled acceptance of arbitrary
    subsets, it is obvious that one should accept choice
    sets.
  \end{quoti}
  Lavine~\cite[page 4]{Lavine:understandinf}:
  \begin{quote}
    The principle [AC] really is inherent in the notion of an
    arbitrary function.
  \end{quote}
  
  \section{Comparisons and challenges}  \label{sect:chall}

  We shall now look at some points of disagreement (or apparent
  disagreement) between my argument and other views
  in the literature.
  
 \subsection{Appeals to intuition}

Firstly, my argument 
 maintains that mathematics is based (at least in part) on intuition; but many authors
 are wary of appeals to intuition.  This wariness sometimes leads them to seek
 other kinds of justification for mathematics.  Some of the resulting projects
 are mentioned in the following sections.

 

\subsection{The need for object realism}

Apart from the role given to intuition, another potentially
controversial aspect of my argument is that it lumps together various kinds of
realism.  As explained
in~\cite{Linnebo:platonismphilmath}, there is a distinction
 between \emph{truth value
  realism}, the belief that sentences of a certain kind are bivalent,
and \emph{object realism}, the belief that ``abstract mathematical
objects exist''.\footnote{Cf.~the view attributed by
  Dummett to Kreisel~\cite{Sundholm:Kreisdict}: ``The problem is not the existence of mathematical objects but the
  objectivity of mathematical sentences.''}
  My narrative maintains that the two go together: for each kind of
  sentence, we should either hold both these beliefs or hold neither.


  For a contrasting view, let us take a look at the ``modal structuralist'' project, which attempts to
justify truth value realism without relying on object
realism.\footnote{A somewhat similar project is presented in Chihara~\cite{Chihara:constructmathexist}.}   In Hellman's account~\cite{Hellman:mathwithout}, the
truth of an 
arithmetical sentence $\phi$ is explained as the  validity of (a
version of) 
the sentence: ``Necessarily, in any second-order PA
model, $\phi$ is true.''

Although this does result in a bivalence principle, I shall raise
two concerns.

Firstly, for a translation between two languages to count as an
explanation, the target language needs to be more meaningful \emph{a
  priori} than the source.  Is this really the case for the modal
structuralist interpretation of PA, whose target language includes
modalities, first-order quantifiers of unspecified range (or perhaps
universal range), and second-order quantifiers?


Secondly, the modal structuralist theory adopts as an axiom the possible 
existence of a second-order PA model.   But how can
this be justified without appealing to the Arbitrary Natural Number
intuition?  

Because of such concerns, I am unpersuaded by this kind of attempt to
sustain 
truth value realism without object
realism.  But let me repeat: construing the latter as something
less than essentialism (Section~\ref{sect:platnotess}) makes it easier
to accept than it would otherwise be.  For related discussion, see
e.g.~\cite{Parsons:maththoughtobj,Linnebo:platonismphilmath,Maddy:defendax,Resnik:mathscipatt,Shapiro:philmathstructont}.

\subsection{Anti-realist theory of meaning} \label{sect:antireal}


We turn next to Dummett's ``anti-realist'' position,\footnote{I present an
  extreme version, for expository convenience. In Dummett's writings,
  many variations are considered.} which concerns not
only mathematics but many other areas of interest.
See~\cite{Loux:realantirealdummchall} for an overview.

According to this philosophy,  
 the bivalence of Googolplex Goldbach should be accepted, but not that of the  
Goldbach conjecture or even the Cleopatra
hypothesis~\cite{Dummett:realpast}.    
These views are based not on a finitistic or constructive ontology (as
 Section~\ref{sect:intu} might suggest), but on a verificationist 
``theory of meaning'' aligned with intuitionistic logic.  An early paper~\cite{Dummett:truth} declares:
\begin{quoti}
  We are entitled to say that a statement $P$ must be either true or
  false [\ldots] only when $P$ is a statement of such a kind that we
  could in a finite time bring ourselves into a position in which we
  were justified either in asserting or denying $P$.
\end{quoti}
This position seems to undermine the premise of our questionnaire.
For it says that, if a sentence
cannot be known ``in a finite time'' to be true or to be false, then
\emph{ipso facto} its bivalence cannot be asserted.  So there is no
need to wonder whether (for example) the totality $\pset\pset\nats$ is
real; the
theory of meaning does all the work.

I find this rather hard to swallow.  Surely questions of mathematical reality    
are genuine and cannot be avoided by mere linguistic
convention?  
In
any case, the anti-realist position and its elaboration by Dummett have given rise to a
substantial literature, both critical and supportive,
e.g.~\cite{Edgington:meanbivreal,McDowell:mathplatdumantireal,Burgess:dummettcaseintu,Wright:realmeantruth,Prawitz:meaningtheoryantireal,Burri:residreal}.  Note, in particular, 
Raatikainen's comparative analysis of intuitionistic notions of
truth and the challenges they face~\cite{Raatikainen:conctruthintu}, 
and Rumfitt's defence of certain kinds of classical
reasoning that Dummett did not accept~\cite{Rumfitt:boundarythought}.

\subsection{The ability to reflect}

In Section~\ref{sect:cons}, I stated that a finitist will accept the consistency of
PRA, and a countabilist that of PA.  This appears to be at odds with
Tait's view~\cite{Tait:finitism} that finitistic reasoning is 
limited to PRA, and with Isaacson's view~\cite{Isaacson:arithmeticalhidden} that
countabilistic (``arithmetical'') reasoning is limited to PA.

The source of the discrepancy is that I have tacitly granted everyone the ability to
reflect on their own language and reasoning.   An \emph{unreflective}
finitist would accept every PRA
proof presented to them in their mother tongue, but (according to
Tait) would consider PRA as a whole to be an unintelligible formal system whose consistency is
doubtful.  It is \emph{unreflective} finitistic reasoning that, in Tait's
view, is limited to PRA.  Likewise it is \emph{unreflective} countabilistic
reasoning that, in Isaacson's view, is limited to PA.  
When the positions of Tait and Isaacson are interpreted in the way
just described, there is no disagreement.

\subsection{Inductive inference}

In Section~\ref{sect:justbelief}, I stated that inductive inference is not accepted,
not even for Liminal Goldbach, let alone for consistency
statements.  Although this is the standard attitude of the
mathematical community, there is a substantial literature that takes
inductive evidence seriously.  For example, such evidence has  
been used in set theory to defend the 
  AD$^{\mathrm{L}(\mathbb{R})}$ 
  hypothesis~\cite{Maddy:believeI,Maddy:believeII,Martin:mathevid},
  in computational complexity theory to defend the $\mathbf{P} \not=
   \mathbf{NP}$ hypothesis~\cite{Aaronson:pnp}, and in many other
   fields of mathematics~\cite{Aberdein:evidproofderiv}.
   
See Paseau~\cite{Paseau:arithenumindsizebias} for a discussion that
specifically considers the inductive evidence for consistency statements.

\section{Conclusions} \label{sect:conc}

We began by taking the view that only proof and intuition can provide
adequate grounds for belief.  This led us from the broad question of what to believe to the more focused question of which intuitions to accept.  

As a result, we have only a few options.  We can doubt the very notion
of human access to platonic reality, but the price of such extreme
scepticism is doubting the bivalence of Googolplex Goldbach---the
ultrafinitist view.  Moving on, we can doubt the bivalence of CH or
the truth of AC, but then the Googolplex consistency of \pathree{} 
is in doubt.  So we cannot adopt a foundational theory
that includes \pathree{}, such as ZF; nor even one that includes
I\pathree{}, such as IZF~\cite{Crosilla:czfizf} or free topos theory.

The focus of this article is higher-order arithmetic, rather than set
theory.  So we have not examined 
particularism and totalism.  Particularists (like me) need to know:
what is the price of doubting the bivalence of GCH?  Totalists
 need to
know: what is the price of doubting the bivalence of the Club-Failure Hypothesis?  These questions are left to future work.

Although the story is unfinished, then, let me sum up.  When setting
out our fundamental mathematical beliefs, we are free to decide how
credulous or sceptical to be.  Our decision will depend on the
strength of the intuitions we experience, and on our degree of
caution.  But this freedom has limits, as we must fully accept or
fully doubt each intuition.  Thus, believing in the consistency of everything and the reality of nothing is not an option.  Scepticism always comes at a price.

\paragraph*{Acknowledgements} I thank the MathOverflow
community for technical help, and Salvatore Florio, Mike Shulman and 
other colleagues and reviewers 
for valuable criticism.

\end{document}